\documentclass[12pt,a4paper]{elsarticle}
\usepackage{amssymb}
\usepackage{amsmath}
\usepackage{amsfonts}
\usepackage{amsbsy}
\usepackage{amscd}
\usepackage{amsthm}
\usepackage{psfrag}
\usepackage{graphics}
\usepackage{color}
\usepackage{colortbl}
\usepackage{graphicx}
\usepackage{subfigure}
\usepackage{array,supertabular}
\usepackage[boxed]{algorithm2e}
\usepackage{tabularx}
\usepackage{booktabs}
\usepackage{multirow}
\usepackage{units}
\usepackage{subeqnarray}
\usepackage{longtable} 
\usepackage{rotating}
\usepackage{subfigure,color,colordvi}
\usepackage[latin1]{inputenc}
\usepackage{url}
\usepackage{latexsym}

\newcommand{\bm}[1]{\text{\boldmath $#1$\unboldmath}}

\newcommand{\domain}{\Omega}
\newcommand{\boundary}{\partial \domain}

\newcommand{\surface}{\bm{S}}

\newcommand{\sref}[1]{Section \ref{section:#1}} 
\newcommand{\fref}[1]{Figure \ref{fig:#1}} 
\newcommand{\tref}[1]{Table \ref{table:#1}} 

\newcounter{myalgorithmctr}

\begin{document}

\begin{frontmatter}
\title{Recent Developments in CAD/analysis Integration}
\author[cardiff]{H. Lian}
\author[cardiff]{S.P.A. Bordas\corref{cor1}}
\ead{stephane.bordas@alum.northwestern.edu}
\author[swansea]{R. Sevilla}
\author[cardiff]{R.N. Simpson}
\cortext[cor1]{Corresponding author}
\address[cardiff]{Institute of Mechanics and Advanced Materials,
  School of Engineering,
  Cardiff University, Queen's Building, The Parade,
  Cardiff CF24 3AA, Wales,UK}
\address[swansea]{Civil and Computational Engineering Centre, College of Engineering, Swansea University, Faraday Building, Singleton Park, Swansea SA2 8PP, Wales, UK }


%
\begin{abstract}
For linear elastic problems, it is well-known that mesh generation dominates the total analysis time. Different types of methods have been proposed to directly or indirectly alleviate this burden associated with mesh generation. We review in this paper a subset of such methods centred on tighter coupling between computer aided design (CAD) and analysis (finite element or boundary element methods). We focus specifically on frameworks which rely on constructing a discretisation directly from the functions used to describe the geometry of the object in CAD. Examples include B-spline subdivision surfaces, isogeometric analysis, NURBS-enhanced FEM and parametric-based implicit boundary definitions. We review recent advances in these methods and compare them to other paradigms which also aim at alleviating the burden of  mesh generation in computational mechanics.
\end{abstract}



\end{frontmatter}
\section{Introduction}
The finite element method is the most widely used numerical method in practice and is underlined by fifty years of rigorous mathematical analysis. It offers results whose quality can be predicted \emph{a priori} and controlled \emph{a posteriori} for a range of problems of engineering relevance.
The finite element method relies on the creation of a data structure known as ``mesh'' which is used to construct the approximation (of both the geometry of the domain and the field variables) and to perform numerical integration. The mesh is a set of elements of polyhedral shape (with straight or curved edges/faces) covering the domain. Traditionally, these elements have been simplex shapes (triangles and tetrahedra) or of quadrangles/hexahedral shapes, but recently developed elements \cite{sukumartabarraei2004} and numerical integration techniques \cite{natarajanbordas2009} allow the development of elements with an arbitrary number of edges/faces.

For day to day industrial problems, initial analyses during early stage design are usually carried out using the assumption of linear elasticity. The geometry of the domain is most commonly generated using CAD software by a designer or engineer. A mesh must then be created to approximate this geometry and discretise the partial differential equation (PDE) governing the problem to allow analysis to be performed. The elements must fulfil a list of quality criteria related to their shape and grading throughout the domain, and this step is still difficult and far from being automated. Once the analysis results are available, it is typically necessary for the analyst to go back to the designer and propose changes to the design so that the part meet certain criteria related to, e.g. the maximum allowable stress or the maximum deflection in the component. In turn, a new mesh must be generated to approximate the new geometry and perform a new analysis. This iterative process must be performed as many times as necessary to converge to a suitable geometry. The idea is illustrated in \fref{iterativeDesign}.

\begin{figure}
\includegraphics[width=0.8\textwidth]{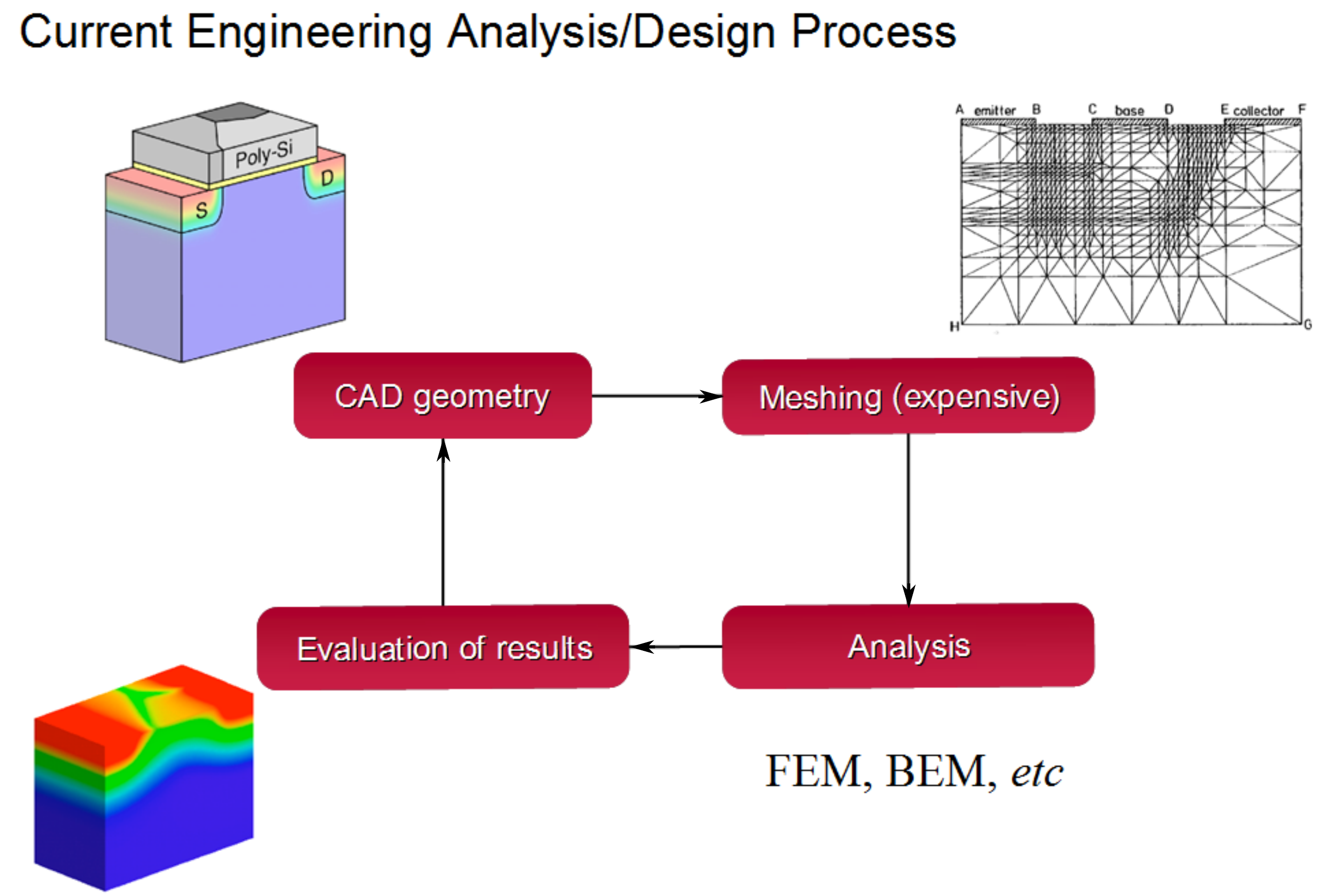}
\label{fig:iterativeDesign}
\centering
\caption[Iterative design process.]{Iterative design process. It is evident that such an iterative process is cumbersome and time consuming because a new, analysis-suitable, mesh must be generated for each new geometry. Several ideas have been advanced in the literature to try and decrease the mesh generation burden. }
\end{figure}

It is evident that such an iterative process is cumbersome and time consuming because a new, analysis-suitable mesh must be generated for each new geometry, which is not a simple task.  Indeed, meshes of good quality elements are still today difficult to generate automatically, or are restricted to simplex (tetrahedral/triangular) elements, which are usually not robust (too stiff and not amenable to incompressible materials). Moreover, it may happen in industrial practice that the geometry is so complex that available mesh generators fail, or require significant human intervention. This has spurred a vast amount of research to enable analysts to generate quality meshes fast and robustly. An exciting avenue of investigation, followed by a large number of strong research groups worldwide consists in developing robust, fast, and automatic mesh generators, in particular for hexahedral elements. Some of the most notable work performed in this area can be non-exhaustively summarised as 

\begin{description}
\item[Triangular/Tethedral-mesh generators] Tetrahedral mesh generators are more flexible and robust for complex geometries compared to hexahedral mesh generators. The current meshing techniques can be roughly grouped into the following categories: advancing front \cite{lohnerparikh1988, loehnerparikh1988}, octree-based methods \cite{yerryshephard1984}, Delaunay approaches \cite{lawson1977, watson1981} and mesh optimization techniques \cite{freitagollivier-gooch1996,cutlerdorsey2004}. Some of the latest developments in this area include works on remeshing approach exploiting the natural anisotropy of most surfaces and  the construction of high-quality isotropic tetrahedral meshes by using a variational principle \cite{alliezcohen-steiner2003, alliezcohen-steiner2005}. The reader is referred to the review paper of Alliez et al. \cite{alliezmeyer2002} for details. The recent paper of Young et al. \cite{youngberesfordwest2008} is particularly noteworthy as it is on the generation of meshes from medical images (image-based mesh generation).
\item[Quadrilateral/Hexahedral-mesh generators] Hexahedral meshes are preferred from a numerical analysis perspective. However, it is more difficult to generate, as explained in the survey of Shepherd and Johnson \cite{shepherdjohnson2008}. Hexahedral meshes can be generated directly, or converted from the tetrahedral meshes indirectly. Zhu \cite{zhuzienkiewicz1991} is one of the first ones using an advancing front approach. Blacker and Stephenson \cite{blackerstephenson1991} introduced a paving algorithm to provide a robust and efficient generator for all-quadrilateral meshes. Plasering is an extension of paving algorithm for hexahedral mesh in three dimensions. The reader is referred to the survey in 1998 by Owen \cite{owen1998} for unstructured meshes and refinement. Recent developments include the work of Zhang and Bajaj \cite{zhangbajaj2006} who proposed an isosurface extraction method to automatically extract adaptive and quality hexahedral meshes from volume data, Ito et al. \cite{itoshih2009} who introduced an octree-based mesh generation method to create geometry-adapted unstructured hexahedral meshes automatically from triangulated surface models, Staten et al. \cite{statenowen2005} who developed an new algorithm to leverage the benefits of paving and plastering, and Sarrate and Huerta \cite{sarratehuerta2011} who proposed an efficient algorithm for complex geometries based on a recursive domain decomposition.
\end{description}

In parallel to these advances in meshing, several paradigms have been advanced in the literature to try and decrease the reliance of numerical analysis on a mesh. The ideas behind these hinge upon the following goals summarised in \tref{meshBurdenIdeas}:
\begin{description}
\item[Avoid altogether the generation of a mesh] Meshfree methods (see Nguyen \cite{nguyenrabczuk2008} for a recent review and computer implementation details including an open source MATLAB code) relax the notion of neighbor relationships (connectivity) by simplifying the addition, suppression and movement of points/nodes. These methods came about in the seventies with the smoothed particle hydrodynamics (SPH) method proposed by Gingold and Monaghan \cite{gingoldmonaghan1977} and  the element-free Galerkin method by Belytschko et al. \cite{belytschkolu1994}. They are not as well understood mathematically as the finite element method, and most methods rely on the definition of a ``mesh'', be it for numerical integration or generation of the point cloud. A rigorous quantification of the approximation property of meshfree methods, especially with reference to the influence of the point distribution is an open problem, although some headway has been made, \emph{inter alia}, in the area of point collocation methods developed by Davydov and Oanh \cite{davydovoanh2011} and maximum entropy interpolants by Rosolen et al. \cite{rosolenmill2010}.

\item[Element technology] Meshes of tetrahedral elements are much more easily generated than meshes of hexahedral elements (bricks) especially if user intervention is to be minimised. If robust and accurate tetrahedral elements could be formulated for most problems of engineering relevance, such fast mesh generators would become all the more useful. Several ideas have been proposed in the literature to avoid most of the limitations of (linear) tetrahedral elements. The node-averaged tetrahedral elements of Bonet and Burton \cite{bonetburton1998}, the locking free tetrahedral element of De Micheli and Mocellin \cite{michelimocellin2009} and the (cell, edge/face, node-based) smoothed finite element method (SFEM) of Liu et al. \cite{liudai2007} are examples of such attempts. Although these methods have other
    drawbacks associated with complications for multi-material problems and larger band widths, the edge-based smoothed finite element method was shown to yield
    superior results to that of the standard FEM for problems in linear elasticity \cite{liudai2007}, visco-elastoplasticity \cite{nguyenliu2009a} and geometrically nonlinear problems \cite{nguyenliu2009}. Polygonal/polyhedral elements are another alternative to ease the difficulties associated with generating quality meshes made up of robust elements. They can serve as transition elements during mesh generation, and therefore offer added flexibility. Work in this area includes the early work of Alwood and Cornes \cite{alwood1969polygonal}, Sukumar and Tabarraei \cite{sukumartabarraei2004}, Euler et al. \cite{euler2006polygonal}, Sukumar and Malsch \cite{sukumarmalsch2006}, Dai et al. \cite{dailiu2007a}, Mousavi et al. \cite{mousavixiao2010} and Natarajan et al. \cite{natarajanbordas2009}.

\item[Use the geometry data provided by CAD directly in analysis] From the description of the early-stage design process above, it is clear that the geometrical information provided by CAD is lost during the mesh generation process, because the smooth, arbitrary surfaces representing the boundary of the object are fitted as closely as possible during the mesh generation by piecewise polynomial surfaces (usually piecewise linear or quadratic) or the faces of the finite elements. The idea of subdivision surfaces B-spline finite elements of Cirak et al. \cite{cirakortiz2000,cirakortiz2001,cirakscott2002} and isogeometric analysis of Hughes \cite{hughescottrell2005} is to use the geometrical description of the component to approximate the field variables. More precisely, the (usually Non-Uniform Rational B-Spline -- NURBS) description of the boundary of the domain is used to build shape functions for the approximation of the unknown fields on the domain. In this way, any modification of the geometry of the object is immediately translated into a modified approximation over the new geometry, without needing to reconstruct a polynomial fit of the boundary. The first set of such methods were developed in the context of the finite element method \cite{hughescottrell2005}, where a domain mesh is required to construct the NURBS-shape functions, and, more recently, for the boundary element method (BEM) in the paper of Politis et al. \cite{politisginnis2009} and Simpson et al. \cite{simpsonbordas2012}, where mesh generation is completely avoided, bearing in mind the intrinsic limitations of the BEM. Advantages of the ``isogeometric approach'' include the use of the \textit{exact geometry at all stages of the analysis} and, more significantly, a significant reduction or circumvention of mesh generation. But recent research also illustrated several other advantages and difficulties of the IGA approach. It is the goal of the present paper to give an overview of some of these developments. Recently, the work of Moumnassi et al. \cite{moumnassibelouettar2011}, based on the initial idea of Mo{\"e}s et al. \cite{moescloirec2003} and Belytschko \cite{belytschkoparimi2003} showed that it was possible, using a marching method, to obtain an implicit domain definition based on multiple level set functions from  arbitrary parametric surfaces provided from CAD data. This was a significant advance compared to previous work, as it enabled the treatment of corners and sharp edges (at the cost of a multiple-level set formulation). In some sense, this class of ambient space methods can be thought of as particular cases of the techniques described in the next paragraph.

\item[Use a structured mesh independent of the geometry] An important difficulty in the mesh generation process emanates from the requirement of the mesh to conform to the (usually complex) geometry of the domain. An alternative is to allow the geometry of the boundary of the domain to be non-conformally represented by the mesh. The boundary cuts the background mesh arbitrarily, and the latter is locally refined, usually using oct-tree based techniques -- which also simplify parallelisation, in particular using modern graphical processing units (GPUs). The method was recently revisited within the context of the extended finite element method \cite{belytschkoblack1999, belytschkoparimi2003, moescloirec2003}.

\item[Use methods only requiring a boundary discretisation] This type of methods take semi-analytical ways to decrease the dimension of the problem. The examples are Boundary Element Method initiated by Rizzo \cite{rizzo1967} and Scaled Boundary Finite Element Method by Song and Wolf \cite{songwolf1997}. Boundary element methods rely on fundamental solutions, which are not available for the non-linear problems and heterogeneous materials. On the other hand, Scale Boundary Finite Element Method can be extended to non-linear materials but not suitable for complex geometries where the scaling centre cannot be easily specified.
\end{description}

\begin{sidewaystable}[htdp]
\caption{Several frameworks have been advanced in the literature to try and decrease the mesh generation burden.}
\begin{center}
\begin{tabular}{p{0.12\textwidth}p{0.12\textwidth}p{0.12\textwidth}p{0.12\textwidth}p{0.12\textwidth}p{0.12\textwidth}}
\raggedright Goals & Meshfree methods & IGA & Element technology &Implicit boundary methods& Boundary discretising methods\\ \hline
\raggedright high order continuity & Yes&  Yes & No & No & No\\ \hline
\raggedright alleviate mesh distortion& Yes & Yes & Yes & No & N/A\\ \hline
\raggedright represent geometries exactly& No & Yes & No &No& No\\ \hline
\end{tabular}
\end{center}
\label{table:meshBurdenIdeas}
\end{sidewaystable}

In this paper, we focus on notions related to isogeometric analysis, in particular to the isogeometric finite element method (IGAFEM), the isogeometric boundary element method (IGABEM), the NURBS-enhanced finite element method (NEFEM) and some of the competing methods, in particular immersed/implicit boundary techniques.

The paper is organised as follows: first, the technology that underpins CAGD and which forms the basis of IGA is outlined; next, the basic concepts of IGA and the fundamental differences to conventional implementations are described; recent developments in IGA are then described, emphasising the beneficial properties over traditional approaches and finally, a comparative review with other methods which can also alleviate the mesh generation is given, followed by the comments on future directions and problems that IGA faces.

\section{Computer Aided Geometric Design Technology}
\label{section:CAGD}

Isogeometric analysis (IGA) \cite{hughescottrell2005} has been successful at bringing two disparate research communities together: that of Computer Aided Geometric Design (CAGD) and numerical analysis. As argued above,  converting CAD models into a form suitable for analysis is considered a laborious task, often dominating the design process. IGA offers a new direction, where, instead, a much more direct route is taken allowing analysis to proceed much more simply than before. This is achieved by using the data provided by CAD models \textit{directly} rather than converting it through a preprocessing routine into a form suitable for analysis. IGA has found that the functions which are used to describe CAD surfaces are often directly amenable for analysis where, for example, properties such as Partition of Unity (PUM), linear independence, non-negativity and weak Kronecker delta property can be proven.

To the authors' knowledge, the first work of integrating CAD and engineering analysis can be traced back to the paper by Kagan et al. \cite{kaganfischer1998}, where B-splines were used as  basis functions to represent both the geometry  and the unknown fields. Following this idea, Cirak et al. \cite{cirakortiz2000} proposed a paradigm for thin-shell analysis, but used subdivision surfaces instead of B-splines. These ideas were formalised and generalised by the acclaimed work of Hughes et al. \cite{hughescottrell2005} based on NURBS, which has been extended to T-splines \cite{bazilevscalo2010} recently. The recent work by Nguyen-Thanh \cite{nguyen-thanhnguyen-xuan2011,nguyen-thanhkiendl2011} and Wang \cite{wangxu2011} has incorporated PHT-splines.

The current section is devoted to providing a non-exhaustive overview of the basic technology provided by the computer aided geometric design community, as this technology is central to isogeometric analysis. We attempt to perform this overview whilst keeping in mind that we wish to use the information provided by CAD for analysis purposes. \tref{comparisonCAGDforAnalysis} gives a comparison of the technology reviewed here, namely B-splines, non-uniform rational B-splines (NURBS), analysis-suitable T-splines, PHT-splines and subdivision surfaces. Other very recent techniques which are not reviewed here include RHT-splines.

\begin{sidewaystable}[htdp]
\caption{Comparison of various CAD descriptions from a geometrical and analysis point of view.}
\setlength{\tabcolsep}{2pt}
\begin{center}
\begin{tabular}{p{0.22\textwidth}p{0.10\textwidth}p{0.10\textwidth}p{0.10\textwidth}p{0.10\textwidth}p{0.10\textwidth}}
Goals & B-splines & NURBS & Analysis-suitable T-splines & PHT-splines & Subdivision surfaces \\ \hline
\raggedright Topology flexibility & No & No & Yes& Yes& Yes\\ \hline
\raggedright Watertight geometry& No & No & Yes& Yes& Yes\\ \hline
\raggedright Smoothness/Continuity & desired continuity & desired continuity&desired continuity&$C^2$& $C^2$ or $C^1$\\ \hline
\raggedright Partition of unity  &Yes & Yes & Yes& Yes& N/A\\ \hline
\raggedright Kronecker delta property & weak & weak & weak & weak & N/A\\ \hline
\raggedright Local refinement & No & No & Yes& Yes&Yes\\ \hline
\raggedright Linear independence & Yes & Yes & Yes& Yes& N/A\\ \hline
\end{tabular}
\end{center}
\label{table:comparisonCAGDforAnalysis}
\end{sidewaystable}%

\subsection{B-Splines}
\label{section:Bsplines}
B-spline geometry is a mapping from parametric space to physical space through a linear combination of B-spline basis functions, which are defined in parametric space, and the corresponding coefficients which are called control points because their physical meaning are a series of points scattered in physical space. B-splines in one dimension can be expressed as:
\begin{equation}
\mathbf{C}(\xi)=\sum_I^{n} N_{I,p}(\xi)\mathbf{B}_I
\end{equation}
where $\mathbf{C}(\xi)$ denotes the physical curve we are interested, $\xi$  the coordinate in parametric space, $B_I$ the control points, $N_{I,p}$ the B-spline basis functions of order $p$.\\
B-splines possess the following properties
\begin{description}
\item[The convex hull property] The B-splines geometry is contained in the convex hull constructed by the control grid, which is a mesh interpolated by control points.
\item[The variation diminishing property] No plane has more intersections with the curve than it has with the control grids. This property renders B-splines less oscillatory than Lagrangian polynomials.
\item[The transformation invariance property] An affine transformation to a B-splines curve can be achieved by applying an affine transformation to the control points.
\item[The partition of unity property]
\begin{equation}
\sum_{I=1}^{n} N_I(\xi) = 1
\end{equation}
\item[The non-negative property] Each basis function is pointwise nonnegative over the domain.
\item[Ease of refinement and control of continuity] It is possible to use the intrinsic building blocks of NURBS for local refinement in the standard $h$ and $p$ way, but also by elevating the order of basis functions and continuity, simultaneously, also known as  $k$ refinement (see \sref{methods-refinement}).
\item[Locally supported] The basis functions are pointwise non-negative in $p+1$ knot spans. This property leads to a fact that the modification of a control point only influences the geometry locally.
\item[Weak Kronecker delta property] A weak Kronecker delta property means $N_I(\mathbf{x})=0$ but $N_I(\mathbf{x}_I)\neq 1$, which is useful for enforcing boundary conditions in engineering analysis, because only the control points corresponding to boundaries need to be considered.
\end{description}

\subsection{Non-Uniform Rational B-Splines}
\label{section:NURBS}

Non-uniform Rational B-Splines (NURBS) \cite{piegltiller1997,scottli2011} are an extension of B-splines by introducing the weights for basis functions to represent a wider range of geometries such as conic sections. From a geometric point of view, NURBS in $\mathbb{R}^d$ is a projective transformation of B-splines in $\mathbb{R}^{d+1}$. NURBS inherit the properties of B-splines as mentioned above, but still have some drawbacks:

\begin{description}
\item[Rational functions] As NURBS are not polynomial functions, integrating them cannot be done exactly using Gau\ss{} quadrature.
\item[Tensor product] The parametric space and control points rely on a structured grid due to the tensor product property of NURBS, which increases the redundancy of the  degrees of freedom.
\item[Continuity] NURBS usually achieve only $C^0$ continuity between the patches.
\item[Geometry repair] From a computational geometry point of view, NURBS based geometry always requires some level of repair due to gaps or overlaps of the various patches making up a geometry.\\
 \end{description}

The following subsections summarise a few of the recent advances to alleviate some of those shortcomings.

\subsection{Subdivision surfaces}
\label{section:subdivision-surfaces}
Subdivision surfaces, firstly proposed by Catmull and Clark \cite{catmullclark1978}, can possess a topology flexibility and ability to construct a watertight geometry. Subdivision surfaces use a recursive rule to interpolate or approximate a smooth surface. The methods can be divided into two classes:
\begin{description}
\item[Interpolating Schemes] When meshes are refined, new vertices are added without changing the vertices in the coarser mesh. Consequently, all of the vertices produced
during subdivision will interpolate the surface. Interpolating schemes for quadrilateral
meshes have been introduced by Kobbelt et al. \cite{kobbeltcampagna1998}, while Dyn \cite{dynlevine1990} and Zorin et al. \cite{zorinschroder1996} described interpolating schemes for triangular meshes. In both cases the limit surfaces are $C^1$ but their curvatures do not exist.
\item[Approximation Schemes] In approximation schemes, in contrast to interpolating schemes, the vertices in the coarse mesh need to be recomputed in refined mesh with adding the new vertices. The examples are the schemes proposed by Catmull and Clark \cite{catmullclark1978} and Doo and Sabin \cite{doosabin1978}, which are the first subdivision surface schemes and use quadrilateral meshes. Loop \cite{loop1987} proposed a scheme based on triangular meshes. Theses approximation schemes can achieved $C^2$ continuous, except of some isolated points which are $C^1$. After parametrisation, subdivision surfaces can be incorporated into IGA \cite{cirakortiz2000,cirakortiz2001,cirakscott2002}.
\end{description}

\subsection{T-Splines}\label{section:t-splines}
T-splines were proposed by Sederberg et al. \cite{sederbergzheng2003} to overcome the drawbacks of NURBS. T-spline control grids, also called T-meshes, permit T-junctions which are similar to the concept of the ``hanging nodes'' and oct/quad-tree meshes in FEM. In this way, lines of control points need not traverse the entire control grid. If a T-mesh is simply a rectangular grid with no T-junctions, the T-spline reduces to a B-spline. T-splines support many valuable operations within a consistent framework, such as local refinement, and the merging of several B-spline surfaces that have different knot vectors into a single gap-free model.
T-splines form a subset of PB-splines, i.e. the basis functions are defined on local knot vectors. However, the local knot vector of T-splines can be inferred from a T-mesh, instead of being defined arbitrarily.  A posteriori error estimation techniques for local $h$-refinement with T-splines have been very recently introduced by Dorfel et al. \cite{dorfeljuttler2010}. Bezilevs et al. \cite{bazilevscalo2010} have tested T-splines on some elementary two-dimensional and three-dimensional fluid and structural analysis problems.\\
The major drawback of T-splines is a significant implementation complexity. In addition, T-splines lose the property of linear independences, which is useful in engineering analysis. But an analysis-suitable T-splines \cite{lizheng2012}, a subset of T-splines, can always satisfy this requirement.

\subsection{PHT-Splines}\label{section:PHTsplines}
PHT-splines (polynomial splines over hierarchical T-meshes) proposed by Deng et al. \cite{dengchen2006} are
constructed on the basis of T-splines and thus inherit their beneficial properties. PHT-splines use a hierarchical T-mesh which has a nested structure and different levels. PHT-splines define a cubic B-spline basis function in a T-mesh, which is level 0, and then modify the basis functions level by level. PHT-splines possess several main merits compared to T-splines:

\begin{enumerate}
\item The basis functions of PHT-splines are polynomial.
\item The local refinement is simpler compared to T-splines.
\item The conversion between NURBS and PHT-splines is very fast.
\end{enumerate}

However, PHT-splines are only $C^1$ continuous, which are not sufficient to represent complex geometries. PHT-splines contain a closed set (linearly independent) of basis functions, which are important for combining them with FEM.

\section{Isogeometric analysis}\label{section:isogeometric}
CAD models provide a complete description of the geometry of a problem, but for analysis to be carried out there is usually a need to carry out some level of preprocessing to arrive at a model suitable for analysis. This is also the case in isogeometric analysis (IGA), although the chief motive of this method is to simplify this preprocessing step.

The difference between tradition practice and IGA however, is that the time and difficulty of this preprocessing step is significantly reduced. For example, in the case of IGA coupled with immersed and boundary element methods, the task of meshing is almost completely circumvented. This section aims to give an overview of the preprocessing steps required for IGA for various numerical methods and the differences and advantages over traditional techniques.

By utilising the exact geometry and noting the advantageous properties of the basis functions that are used in CAD (see \sref{NURBS}), several important developments have been made in the computational mechanics community. Moreover, since the modifications required to convert a traditional numerical analysis code into one that utilises the same basis functions as CAD are often minimal, the method can be adopted quickly and applied to existing numerical methods. Several of these recent developments are highlighted in the present section where the benefits and drawbacks compared to traditional approaches are demonstrated.\\

\subsection{Motivation}

As discussed briefly, the main purpose of isogeometric analysis (IGA), presented by Hughes et al. \cite{hughescottrell2005} is to increase the ties between Computer Aided Design (CAD) data and the numerical method of choice for analysis. In the original paper, this method was the Galerkin Finite Element Method (FEM).

For simplicity, let us take the example of linear elastic structural
mechanics problems. In engineering practice, the geometry is described by
 CAD packages allowing the construction of very complex
three-dimensional (3D) geometries, usually described by several NURBS (Non-Uniform
Rational B-Splines) patches. This geometry then needs to be meshed\footnote{This meshing
operation is inexact, and only as the mesh size goes to zero is the geometry
exactly represented. This geometry approximation error however converges to zero faster than the approximation error on the unknown fields.} and analysed
(e.g. by the FEM) to obtain an approximate distribution of the unknown fields within the
structure. This analysis is used by the engineer to assess the design. If
alterations to this design are required, they must then be communicated to the
drafting team, and the geometry of the component suitably modified. A new mesh
must then be regenerated by the analysis team to perform a second stress
analysis, and the process is repeated as many times as is required for a suitable
design to be obtained.

Considering this design cycle, it becomes clear
that tying the geometrical CAD data to the mesh in an automatic and direct way,
where the mesh would be defined automatically from the CAD data is highly
desirable. This desirable integration is evident from the consolidation trends
observed in the CAD/Analysis industry where CAD companies join forces with
leaders in analysis, with the goal to streamline the design cycle and decrease
lead time in computer aided design.

In IGA the geometry can be
represented \emph{exactly} for analysis by direct use of CAD data through the use of
NURBS-based approximations. This approach is being pursued widely both in the
engineering community by improving the basis functions and seeking adaptive
methods \cite{liptonevans2010,bazilevscalo2010,dorfeljuttler2010} and applying the
ideas to fields outwith mechanics, such as electromagnetics
\cite{buffasangalli2009}. The idea is also being taken up by the applied mathematics
community through the derivation of error estimates by, e.g. \cite{beiraodaveigabuffa2011b}. This
engineering and mathematics research has allowed a rather rapid development of the method
through application-oriented investigations and sound mathematical analysis,
respectively.

\subsection{Methods of refinement}
\label{section:methods-refinement}
One refinement method is known as knot insertion, which is analogous to $h$-refinement in FEM and relies on
splitting an element. However, the insertion of the existing values can reduce the continuity
of basis function without producing new elements. The second is order elevation, which has similarities
with $p$-refinement in classical FEM. Finally, $k$-refinement has no analogue in FEM, because it allows
increasing the order of the basis functions and the number of elements, simultaneously. The advantage of this method
is that fewer degrees of freedom are introduced compared to $p$-refinement. \cite{cottrellhughes2007}
details the $h$-$p$-$k$- refinement schemes, Dorfel et al. \cite{dorfeljuttler2010} explains an adaptive $h$-refinement, and Beir{\~a}o da Veiga et al. \cite{beiraodaveigabuffa2009} introduce error estimates for $h$-$p$-$k$- refinement schemes.

\subsection{A priori error estimation}
\label{section:errorEstimation}

The difficulties of a priori error estimation in isogeometric analysis arises from two aspects:
\begin{enumerate}
\item NURBS are not interpolatory.
\item The support of NURBS basis functions are in general found to be larger than polynomial approximations of an equivalent order.
\end{enumerate}
 The first difficulty is solved by pulling back the approximated function to the regular parametric space and noting that NURBS in $\mathbb{R}^d$ is a projective transformation of B-splines in $\mathbb{R}^{d+1}$. The second difficulty is addressed by introducing the so-called ``bent'' Sobolev spaces \cite{bazilevsdaveiga2006}. A meaningful approximation expression is derived to include not only the approximated function but also the gradient of the mapping, which does not change with the refinement.

\emph{Finally, a significant result is that: the IGA result using $p$-order NURBS has the same order of convergence  as the Finite element method using $p$-order polynomials.}

\subsection{B\'{e}zier extraction}
\label{section:bezier-extraction}
To further integrate isogeometric analysis with existing FEM codes, B\'{e}zier extraction is utilised for NURBS \cite{bordenscott2011} and T-splines
\cite{scottborden2011} to allow numerical integration of smooth functions to be performed on $C^0$ B\'{e}zier elements. The idea of B\'{e}zier extraction is that localised NURBS or T-spline basis functions can be represented by a linear combination of the Bernstein polynomials. In addition to localising the support of basis functions into an element, B\'{e}zier extraction provides an
element data structure suitable for analysis. That is,  similar to Lagrangian polynomials in traditional FEM, Bernstein polynomials do no change element by element. Furthermore, no intermediate
parametric space is employed, hence the physical geometry can be mapped directly from parent elements. However, it should be noted that B\'{e}zier extraction increases
the degrees of freedom of the system.

\subsection{Isogeometric boundary element method}
\label{section:boundElemMethod}

As mentioned previously, the main bottleneck of IGA is the transition from the boundary representation to a domain representation suitable for analysis. An isogeometric approach using the framework of the boundary element method - coined the isogeometric boundary element method (IGABEM) - was proposed to solve this problem. The idea relies on the fact that both CAD models and boundary element methods rely on quantities defined entirely on the boundary. In the IGABEM, NURBS (and later, T-splines) basis functions are used to discretise the boundary integral equation (BIE) for the surface geometry, boundary displacement and traction fields. The first implementations of this concept were published by Simpson et al. \cite{simpsonbordas2012} and Politis et al. \cite{politisginnis2009} for two dimensional linear elastostatic problems and exterior potential flow problems, respectively.

 Recently, IGABEM was extended to include T-splines for 3D linear elastostatic analysis in \cite{scottsimpson2012} with analysis on complex geometries demonstrated. Moreover, Scott et al. \cite{scottsimpson2012} take advantage of a regularised form of BIE to simplify the implementation of singular integration inherent in boundary element methods. IGABEM inherits the benefits of both IGA and the classical boundary element method representing an attractive technique for early-stage design.

 IGABEM also promises important advantages for problems with infinite domains, such as acoustic and electromagnetic problems. Advantages are also anticipated in areas such as moving boundary problems (including fracture and contact) and inverse problems. Similar to classical BEM, however, IGABEM is not naturally adapted to non-linear problems such as elasto-plasticity because of the requirement to evaluate a domain integral and thus, the application of the method is mainly restricted to linear problems. Compared to IGAFEM, the mathematical foundations of IGABEM require further investigation, particular with reference to methods that utilise collocation.

\subsection{Analysis-aware geometric technologies}

IGA circumvents some of the difficulties associated with  meshing but still have to overcome severe problems associated with CAD model quality. A main challenge of IGA is that most of the currently available methods, to the exception of the isogeometric boundary element method \cite{simpsonbordas2012} described in \sref{boundElemMethod}, are based on a volumetric model, whilst CAD models only provide a surface model. Furthermore, commercial CAD systems do not provide boundary representations of models in any general way to help construct full volume representations. The work to generate a trivariate geomertical representation is necessary. Interested readers can refer to \cite{cohenmartin2010,lu2009,kimseo2009,costantinimanni2010,aignerheinrich2009,martincohen2009,wangzhang2011,xumourrain2011,xumourrain2011a}.\\

\subsection{Applications}
\subsubsection{Plates and shells}
There are two kinds of theories for plate and shell problems. One is the Kirchhoff-Love (KL) theory, which is suitable for thin plates or thin shells. In Kirchhoff-Love theory, rotations, which are the first derivatives of the displacements, are not interpolated directly, but through the displacement approximation, hence leading to the requirement of a $C^1$ continuous representation approximation of the displacements, to yield $C^0$ rotations. To construct a $C^1$ element is not easy for polynomial basis functions. The second is Mindlin-Reissner theory for thick plates and shell problems, where shear-locking has to be dealt with carefully. In addition, the number of degrees of freedom is increased compared to the KL theory because rotations need additional interpolation.

A rotation free analysis can be performed by IGA, i.e. the high order continuity between elements can be achieved by only interpolating the displacement without approximating its derivatives. Another apparent advantage of IGA for shells is that it can capture the exact geometry to enhance the analysis accuracy. Although IGA for thin plates (shells) require a rotation boundary condition and continuity at finite angles,  these constraints can be enforced by Lagrange multipliers or other similar methods. The application of IGA to plates and shells can be found in \cite{bensonbazilevs2010, kiendlbletzinger2009, bensonbazilevs2010a, beiraodaveigabuffa2011, uhmyoun2009}, to which the interested reader is referred to for details.
\subsubsection{Fracture and damage}
\label{section:fracture}
The traditional Finite Element Methods have difficulties dealing with crack propagation problems because of the mesh regeneration adapted to the moving boundaries and the singularity near the crack tip. A great success in this area has been achieved by extended finite element methods proposed by Belytschko and Black \cite{belytschkoblack1999}. The property of partition of unity of NURBS basis functions can lead to a combination of IGA and the extended finite element method. This combination is especially potent for complex geometries and curved cracks. First results were presented in the papers by De Luycker et al. \cite{luyckerbenson2011} and Ghorashi et al. \cite{ghorashivalizadeh2012}. Moreover, IGA also possesses advantages in implicit gradient damage models which necessitates a higher order fields. This problem has been circumvented by the use of meshless methods as in, e.g. \cite{askesaifantis2001,paminaskes2003,borst2009}. Isogeometric finite elements allow for the construction of higher-order continuous field on complex domains, which was first demonstrated by Verhoosel et al. \cite{verhooselscott2011a}. Other advances of the application of IGA in this area were presented by Borden et al. \cite{bordenverhoosel2012} using a phase field model, de Borst et al. for failure \cite{borsthughes2011} and Verhoosel et al. \cite{verhooselscott2011} for cohesive fracture.

\subsubsection{CFD and CEM}
%
The importance of the geometrical model in the numerical solution of Euler equations of gas dynamics in the presence of curved walls was identified in a Finite Volume context by~\cite{Dadone-DG:94,Barth:98}. The same problem when solving Euler equations was identified in a discontinuous Galerkin framework in~\cite{Bassi-BR:97} and a detailed study presented in~\cite{CMAME-Ven-02} concluded that accurate results can only be obtained taking  into account the curvature of the boundary. The benefit of using an exact boundary representation in the numerical solution of compressible flow problems was illustrated in~\cite{IJNMF-NEFEM,ARCME-NEFEM}, showing that the error in aerodynamic quantities of interest can be reduced significantly when the approximated boundary representation of standard isoparametric finite elements is replaced by the NURBS boundary representation.
The effect of accurate boundary representations in the numerical solution of viscous problems has been also demonstrated by~\cite{DGcurvedViscousBook,DumbserTurbulentCurvedBoundary}, but the importance of the geometrical model is not exclusive of fluid mechanics. Maxwell's equations are also very sensitive to an accurate geometric description. In~\cite{xuedemkowicz2005} the error induced by isoparametric approximations of curvilinear geometries is studied. By solving the 3D Maxwell's equations in a sphere, the authors show that an exact description of the geometry reduces the error in one order of magnitude compared to isoparametric elements. In~\cite{sevillafern'andez-m'endez2008,sevillafern'andez-m'endez2011a,sevillafern'andez-m'endez2011} several electromagnetic scattering applications are used in order to demonstrate the benefits of using the exact boundary representation of the domain instead of the classical isoparametric description. In every example shown, IGA is  more accurate than the standard isoparametric approach because of the geometrical accuracy. The application of IGA in this area can refer to the work of \cite{akkermanbazilevs2008,bazilevscalo2007, bazilevsakkerman2010, bazilevsmichler2010} in fluids and turbulence, \cite{bazilevscalo2008,gomezhughes2010,nielsengersborg2011,bazilevscalo2006,zhangbazilevs2007,bazilevsgohean2009} in fluid structure interaction, and \cite{buffasangalli2009,vazquezbuffa2010} in electromagnetics.

\subsubsection{Contact}
\label{section:contact}
Contact problems are still a challenging topic in computational mechanics because of the nonlinear problems arising from the boundary and the high sensitivity to the geometry accuracy, especially when sliding is involved. IGA satisfies the geometric accuracy requirement, however, the schemes for enforcing the contact boundary condition must be chosen carefully due to the fact that the control points are not located on the boundary. The methods fall into two classes. The first is the non-mortar method, which exerts contact boundary condition at the quadrature point. This scheme is also called knot-to-surface method in IGA because each quadrature point is represented by a parametric coordinate. On the other hand, mortar methods are a segment-to-segment technique and the boundary condition is imposed weakly with Lagrangian multipliers. The first work dealing with contact treatment in IGA includes \cite{wriggerszavarise2007,temizerwriggers2011,delorenzistemizer2011,kimyoun2011,lu2011, delorenziswriggers2011}.

\subsubsection{Waves and vibrations}
\label{section:vibrations}

One of the advantages of the IGA for vibration problems is that the NURBS basis functions are non-negative, thus the mass matrix is also non-negative. In addition, the finite element method results in a so-called acoustical branch and an optical branch.  This branching is due to the fact that there are two distinct types of difference equations for the finite element: those corresponding to the end-point nodes at element boundaries, and those corresponding to mid-point nodes on element interiors. The acoustical branch corresponds to modes in which the neighboring end- and mid-point nodes oscillate in phase with each other, and the optical branch modes are the modes in which they are out of phase. Alternatively, the quadratic NURBS difference equations are all identical, and no such branching takes place. The reader is referred to \cite{cottrellreali2006} and \cite{hughesreali2009} for details.

\subsubsection{Structural shape and topology optimisation}\label{section:optimisation}

In a sense, IGA seems very suitable for the shape optimisation because of the tight interaction with CAD. Indeed, whether the optimisation is done manually, by trial and error, or automatically, through optimisation algorithms, the fact that any change in the shape automatically results in a change in the approximation of the field variables, without any additional mesh generation, makes the method attractive. Furthermore, the design velocity is accurate and the result does not require any smoothing process. Three of the relevant very recent papers in this area are \cite{wallfrenzel2008,choha2009,seokim2010}. Nagy et al. \cite{nagyabdalla2010} study the variational formulation of stress constraints in isogeometric design.

A drawback of this method, which is a direct byproduct of its advantages, is that the geometry approximation is intrinsically tied to the field approximation. In such a formulation, the flexibility of modifying the field approximation independently of the geometry approximation is lost. Although no detailed study of the interrelationships between geometry error and field approximation error during structural optimisation is known to the authors, this flexibility may be important to ensure convergence to a globally optimal solution since it is not clear from the outset that the geometrical basis is optimal from an analysis standpoint.

For topology optimisation, the initial method is cell-based or element-based, where the geometry can be represented by level sets \cite{oshersethian1988} or splines. The traditional FEM encounters several main challenges. Firstly, in the fixed grid method, the optimal result depends severely on the design space which is constructed by the fixed grid. Secondly, postprocessing is required to transfer the end result to a CAD-suitable model, unless NURBS or other CAD-suitable approximations are used. Recent work in this direction was presented by \cite{bugedarodenas2008}. NURBS-based IGA can avoid the difficulties of cell-based methods by using splines to represent the geometry. However, a trimmed-NURBS technique is necessary to deal with the topology changes \cite{seokim2010, seokim2010a}.
\subsubsection{Locking}
It is well known that depending on the approximation properties of the numerical method, volumetric locking arises when dealing with (nearly) incompressible materials. A simple and effective method to overcome this difficulty is to use standard Gau\ss{} integration for the volumetric strain and a lower integration order for the shear strain. However, this method is not suitable for anisotropic and large deformation problems.

A more general way is $\bar{\mathbf{B}}$ method proposed by Hughes \cite{hughes1980}, which is a three field method where the strain tensor is split into a deviatoric and dilatational parts, and the latter part is approximated by its linear projection into lower order space. Elguedj et al. \cite{elguedjbazilevs2008} have investigated $\bar{\mathbf{B}}$ in isogeometric analysis for linear incompressible elasticity and $\bar{\mathbf{F}}$, a extension of $\bar{\mathbf{B}}$, for nonlinear elasticity, and concluded that IGA is more robust than the traditional $p$-methods.

Concerning shear locking, NURBS finite elements are not necessarily locking-free but this can be alleviated by a discrete shear gap method, as described in \cite{echterbischoff2010}.

\subsubsection{Phase field and gradient theory}
Both phase-field methodology and gradient theory typically lead to high order partial-differential equations. The Cahn-Hilliard equation governing phase field problems is a fourth-order differential equation. Consequently, the weak form necessitates $C^1$ continuous functions. To represent the size effect which cannot be explained by classical elasticity theory, the gradient theory was first proposed by Mindlin \cite{mindlin1964,mindlineshel1968}. Following this work, various versions of the method have been proposed but they share the common features that the strain gradient is introduced to the constitutive equation, yielding a higher-order equation than the classical elasticity theory. Traditional methods experience difficulties in constructing high order fields, which can be solved by IGA readily. The application of IGA in phase-field and gradient elasticity problems can be found in \cite{g'omezcalo2008} and \cite{fischerklassen2011} respectively.

\subsection{Challenges}

We have seen in the preceding sections that isogeometric analysis offers distinct advantages over the standard CAD-Analysis-CAD design cycle. The main challenge of isogeometric analysis is to generate a three dimensional discretisation for the interior of the domain when knowing only the boundary representation, given by the CAD model. This volumetric parameterisation has been the topic of a number of recent research papers, and some of them are produced under Marie Curie Initial Training Networks in Europe. The reader is referred to the following work \cite{aignerheinrich2009,martincohen2010}.

In addition to the ``mesh''  generation difficulty mentioned in the previous paragraph, there are a few more stumbling blocks, which require further investigation\footnote{The reader is referred to \tref{meshBurdenIdeas} for a summary. }:

\begin{description}
\item[Integration] Efficient numerical integration of the underlying approximation in isogeometric analysis is difficult since most of those are non-polynomial. This difficulty is common to enriched FEM and most meshless methods. The most recent advance known to the authors is the work of Hughes et al. \cite{hughesreali2010}, where the number of optimal quadrature points is chosen to be independent of the order of basis functions, and roughly equal to half the number of degrees of freedom.
\item[Mesh distortion] The behaviour of isogeometric approximations for distorted meshes has been studied very recently by Lipton \cite{liptonevans2010}, where it is pointed out that current finite element mesh quality metrics are too restrictive for IGA, i.e. that IGA is less sensitive to mesh distortion. Echter and Bischoff \cite{echterbischoff2010} found that higher order continuity between individual patches can significantly improve the accuracy.
\item[Boundary conditions] The lack of strong Kronecker delta property of isogeometric approximations (similar to many meshfree methods) leads to complications in the imposition of essential boundary conditions. A recent paper on this topic is \cite{bazilevsmichler2010}. Note that isogeometric approximations possess a ``weak Kronecker Delta property." This property endows IGA an advantage over meshfree methods in the enforcement of boundary conditions.

\item[Error estimation and adaptivity] As for any numerical methods, a priori and a posteriori error estimates and adaptivity are of prime importance for reliable isogeometric analysis. A method based on $h-$
refinement was proposed very recently by Dorfel et al. \cite{dorfeljuttler2010}. However, $p-$
and $k-$ (see \cite{hughescottrell2005} for a precise definition) refinements are also
possible, which offer great flexibility in the way the approximation space
can be tailored to the unknown solution. One particularly interesting point is the measure of the adequacy of CAD-based approximations for analysis, in particular in the case of topology and structural optimisation.
\end{description}

\section{Other methods to suppress the mesh burden}
\label{section:other-cad-integrated}
\subsection{Meshfree Methods}
\label{section:Meshfree}
It was quickly realised that the finite element method was not best suited to all problems, in particular for large deformation problems and fragmentation. The first meshfree method developed was the smoothed particle hydrodynamics (SPH) by Gingold and Monaghan \cite{gingoldmonaghan1977} in the context of astrophysics problems.

Meshfree/meshless methods also have advantages in the context of CAD to analysis transition, since they allow placing points/nodes with a lot more freedom than in the standard finite element method. This is because, in contrast to finite element methods, they do not employ elements in the construction of the approximation. Instead, a set of nodes, or centres, associated with a domain of influence are sufficient. The connectivity between the nodes, determined by the overlapping of these domains of influence can be defined more flexibly than in the FEM, and, in particular, can evolve in time.

The basic idea behind most meshfree methods is thus to lift the strict connectivity
requirements posed by the FEM. However, the arbitrariness in the node placement is relative since the quality of the
approximation is known to be dependent on the geometrical location of the nodes. A recent and very interesting paper on this matter is that by Davydov and Oanh \cite{davydovoanh2011}. Based on these nodes and the
size of their domain of influence, an approximation is constructed and can be enriched in various ways, to better represent the solution. The reader is encouraged to refer to the recent review \cite{nguyenrabczuk2008} for details and an open-source Code.\footnote{The code, along with enriched isogeometric finite element and boundary element codes may be downloaded here : \url{http://cmechanicsos.users.sourceforge.net/}}.

Some of the most salient points to mention about meshfree methods may be briefly and non-exhaustively summarised as follows:

\begin{itemize}
\item Mesh generation is vastly simplified, although work on the approximation properties of given point constellations is still ongoing. However, the influence of the node placement is generally not as severe as in the FEM.
\item Arbitrary continuity and completeness can be achieved which can be useful for higher order continua or Kirchhoff-Love  plate/shell formulations.
\item Discontinuities can be treated by modifying the weight functions, or through enrichment, similarly to partition of unity enrichment. See the paper of Nguyen et al. \cite{nguyenrabczuk2008} for details.
\item Nodes can be added easily to locally adapt the approximation, as in the work of, e.g. Rabczuk and Belytschko \cite{rabczukbelytschko2005}.
\end{itemize}

Most meshfree methods may be recast as a particular case of the weighted
residual method, with particular choices of the trial and test
spaces. Depending on these choices and whether local or global weak forms are
employed, a variety of methods can be constructed, with their respective
advantages and drawbacks in terms of stability, accuracy, computational cost,
convergence and robustness. This is explained in detail in the papers by Atluri's group, see, e.g. \cite{atlurizhu1998}. \\

IGA and meshfree methods share quite a number of similarities, in particular the advantages of providing high order continuity. Compared to IGA, some meshfree methods (the reader can refer, for more details, to the well-known articles by Belytschko et al. \cite{belytschkoguo2000}, Rabczuk and Belytschko \cite{rabczukbelytschko2004}, and Xiao and Belytschko \cite{xiaobelytschko2005} for a discussion on stability and accuracy of meshfree methods based on Eulerian and Lagrangian kernels) can be considered to somewhat alleviate the difficulties associated with mesh distortion.
On the other hand, IGA retains the exact geometry, and can more easily be incorporated into existing FEM analysis software, as was shown by recent releases of open-source implementations by Vuong et al. \cite{vuongheinrich2010}, de falco et al. \cite{defalcoreali2011} and Nguyen-Vinh et al. \cite{nguyensimpson2012} \footnote{See \url{http://cmechanicsos.users.sourceforge.net/} to download the 3D open source MATLAB code for enriched isogeometric analysis.}. A combination with IGA and meshfree methods is possible by taking advantage of point based splines (PB-splines), as mentioned in \cite{bazilevscalo2010}. 

\subsection{NURBS-Enhanced Finite Element Method (NEFEM)}
\label{section:nurbs-enhanced-fem}

In isogeometric analysis the whole domain is represented as a
NURBS entity, and the solution of the boundary value problem is
approximated with the same NURBS basis used for the description of
the geometry.
There are two main differences between isogeometric analysis and
NEFEM. First, NEFEM considers the exact NURBS description for the
boundary of the computational domain, the usual information provided
by a CAD software. Secondly, NEFEM approximates the solution with a
standard piecewise polynomial interpolation. Moreover, every
\emph{interior} element (i.e.\ elements not having an edge or face
in contact with the NURBS boundary) can be defined and treated as a
standard FE. Therefore, in the vast majority of the domain,
interpolation and numerical integration are standard, preserving the
computational efficiency of the classical FEM. Specific numerical
strategies for the interpolation and the numerical integration are
needed only for those elements affected by the NURBS boundary
representation.

NEFEM was first presented for 2D domains in \cite{sevillafern'andez-m'endez2008},
showing the advantages in front of classical isoparametric FEs using
both continuous and discontinuous Galerkin formulations for the
numerical solution of some test problems. It is important to remark
that all the ideas presented in  \cite{sevillafern'andez-m'endez2008} are valid not
only when the boundary of the domain is parameterized by NURBS, but
for any piecewise boundary parametrization. The discussion is
centered on NURBS boundary parametrization because they are the most
extended technology in CAD. In \cite{sevillafern'andez-m'endez2008} NEFEM was shown to
be a powerful tool for solving compressible flow problems governed
by the Euler equations of gas dynamics.
Several high-order FE methodologies for the treatment of  curved
boundaries are discussed and compared in \cite{sevillafern'andez-m'endez2011a},
including isoparametric FEM, Cartesian FEM, $p$-FEM and NEFEM.
Numerical examples show that NEFEM is not only more accurate than FE
methods with an approximate boundary representation, but also
outperforms $p$-FEM with an exact boundary representation, showing
the advantages of combining Cartesian approximations with an exact
boundary representation. In \cite{sevillafern'andez-m'endez2011} the extension of NEFEM to 3D domains  is
presented. We brief NEFEM in the following paragraph.

Consider an open bounded domain $\domain$ whose boundary
$\boundary$, or a portion of it, is parameterized by NURBS curves in
2D or surfaces in 3D. In 3D every NURBS is assumed to be
parameterized by
\begin{equation*}
 \surface : [0,1]^2 \longrightarrow \surface([0,1]^2)
          \subseteq \boundary \subset \mathbb{R}^3.
\end{equation*}

A regular partition of the domain
$\overline{\domain}=\bigcup_e\overline{\domain}_e$ in elements is
assumed, such that $\domain_i\bigcap\domain_j=\emptyset$, for $i\neq
j$. As usual in FE mesh generation codes, it is assumed that every
curved boundary face belongs to a unique NURBS. That is, one element
face can not be defined by portions of two, or more, different
NURBS. Note however that the piecewise definition of NURBS is
independent of the mesh discretization. Thus, the NURBS
parametrization can change its definition inside one face, that is,
FE edges do not need to belong to knot lines. This is a major
advantage with respect to the isogeometric analysis. Moreover,
different faces of an element can be defined by different NURBS
surfaces. \fref{NEFEM3D_almondSurfaceTriWithPatches} shows
the knot lines of the surfaces that define the NASA almond and the
surface triangulation. It can be observed that the spatial
discretisation is independent of the piecewise NURBS surface
parametrization.

\begin{figure}
\centering
\subfigure[]{\includegraphics[width=0.31\textwidth]{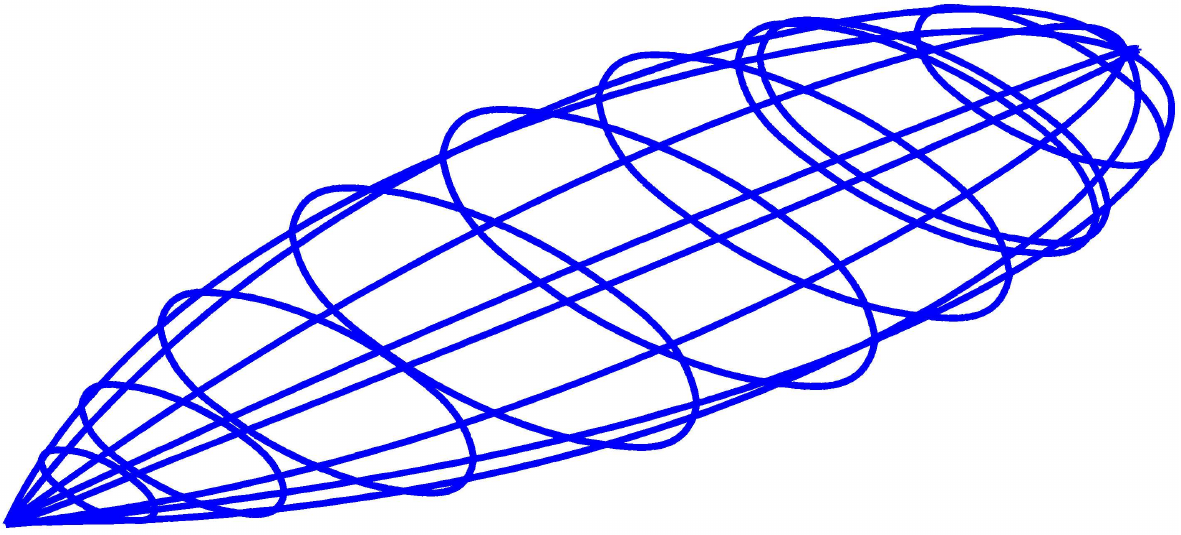}}
\subfigure[]{\includegraphics[width=0.31\textwidth]{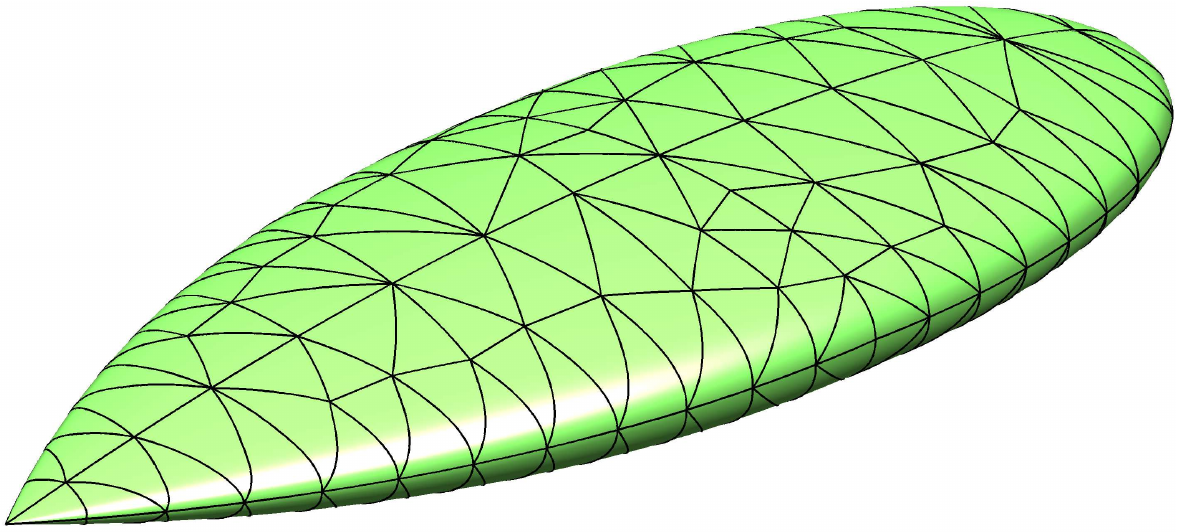}}
\subfigure[]{\includegraphics[width=0.31\textwidth]{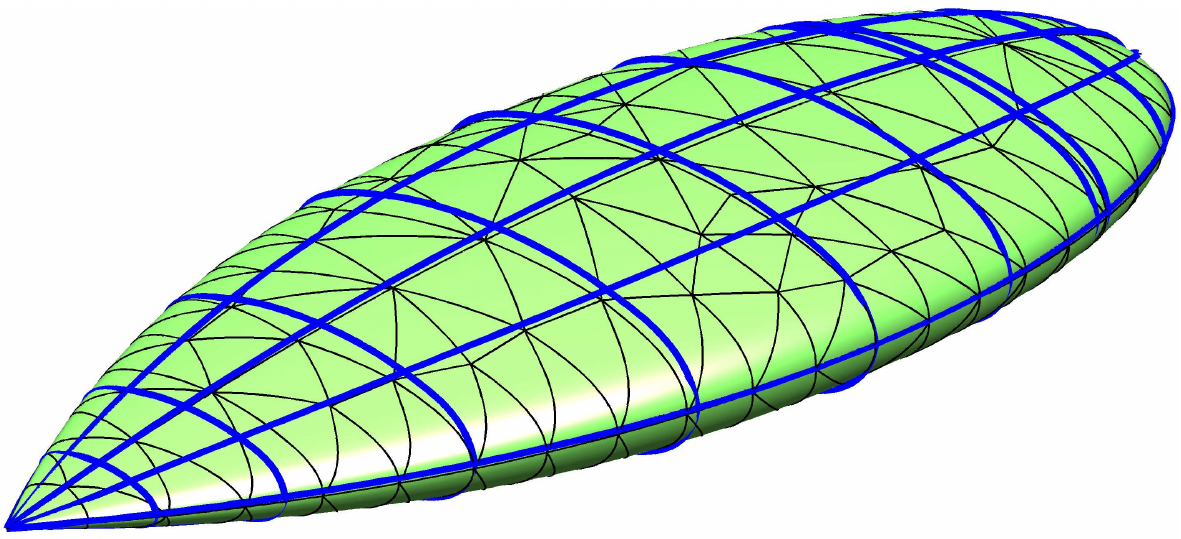}}
\caption{(a) Knot lines of the NURBS surfaces defining the NASA
almond, (b) surface triangulation, and (c) surface triangulation and
knot lines} \label{fig:NEFEM3D_almondSurfaceTriWithPatches}
\end{figure}


NEFEM possess the advantage of accurate representation of geometry and alleviate the difficulty of generating the interior isogeometric element. However, the difficulty of constructing high order continuous field still exist because of the use of polynomials. In addition, the control points are not related to the nodal parameters, hence it is not a geometry ``interaction'' language and loses the advantages in shape optimisation.

\subsection{Implicit Boundary Methods}
\label{section:immersedBMethods}

The implicit boundary method is a class of methods, including immersed boundary methods \cite{mittaliaccarino2005}, fictitious domain \cite{glowinskipan1994}, embedded boundary \cite{johansencolella1998}, virtual boundary \cite{saikibiringen1996} and Cartesian grid methods \cite{yemittal1999}. These methods share the common features, that is, the meshes do not need to conform to the geometry. In addition, its advantage over IGA is that it allows a certain flexibility in the choice of basis functions, which may be different for the field variables and the geometry of the domain. However, it necessitates a ``mesh geometry interaction'' step in order to find the elements which are ``touched'' by the implicit boundary. The differences between various versions of implicit methods depend on how to construct the domain using a background mesh, explicitly or implicitly, the choice of integration scheme and the imposition of essential boundary conditions. Normally, the following three methods are used: Lagrange multipliers, penalty method, Nitsche's method (related to the family of augmented Lagrangian methods) or modifying the basis functions to satisfy the constraint directly.

Moumnassi et al. \cite{moumnassibelouettar2011} propose an algorithm to define the implicit boundary of a domain  from arbitrary parametric representations. More specifically,  an automatic conversion from parametric surfaces to zero level sets for general unstructured meshes is presented by the authors. The use of parametric information can control geometrical errors at the boundaries which affect the convergence rate and solution accuracy, without changing the underlying fixed mesh for analysis. This algorithm helps  integration with computer aided design and manufacturing systems.

\subsection{Smoothed Finite Element Methods}\label{section:sfem}

Recently, novel finite element methods were born from the coupling of stabilized
conforming nodal integration with the standard finite element method
\cite{liunguyen2007}. An overarching theory has been developed in the recent
paper \cite{liu2008}. The main premise of this theory is the wish to achieve
reliable results using lower order elements, i.e. simple meshes (triangles, tetrahedra). SFEM retains the accuracy and inherit the advantages of  triangular and tetrahedral meshes to represent complex geometries and can benefit directly from any advance in automatic remeshing.

 Furthermore, smoothed FEMs are a lot less sensitive to locking (volumetric and shear) as well as mesh distortion (because Jacobians are not required since no isoparametric mapping is used. In this sense, SFEMs are a way to improve the quality of the results obtained by simplex elements, thereby significantly reducing the need for human-intervention in the generation of hexahedral meshes.
%

A few salient features of the smoothed finite element method can be summarised as follows:

\begin{itemize}
\item Reduced sensitivity to mesh distortion.
\item Softer response compared to linear simplex elements (triangular and tetrahedral) thereby allowing the use of unstructured meshes, much easier to generate automatically than structured hexahedral meshes.
\item Possibility to obtain upper bounds of the solution both compared to the displacement-based FEM solution and to the exact solution.
\item A range of methods producing results comprised between quasi-equilibrium methods and displacement-based finite elements.
\item Possibility to obtain ``ultra accurate solutions'' with the ``$\alpha$'' FEM (see e.g. \cite{liunguyenxuan2009} and \cite{nguyenthanhrabczuk2010}).
\end{itemize}

The basic ideas of strain smoothing in FEM and XFEM have been reviewed recently
in \cite{bordasrabczuk2008} and the reader is referred to this paper for more
detail.

\subsection{Solution Structure Method}
Solution structure methods (SSM) \cite{shapirotsukanov2011} can be viewed as a subset of the implicit boundary methods, but the idea to enforce boundary conditions is different and is worth briefing in this section. The key idea of SSM is to use a Kantorovich method, i.e. the solution of a differential equation with homogeneous Dirichlet boundary conditions can be represented in the following form:
\begin{equation}
u=\omega\Phi
\end{equation}
\noindent where $u$ is the solution of the differential equation in which we are interested, $\omega$ is a known function that is positive in the domain and vanishes on the boundary,
and $\Phi$ is some unknown function, which can be discretised by
\begin{equation}
\Phi=\sum_{I=1}^n N_I\phi_I
\end{equation}
\noindent where $N_I$ denotes the basis functions, which can be polynomials, B-splines, radial basis functions. In SSM, $\omega$ is taken as the Euclidean distance function and the basis functions are B-slines. SSM can be thought of as a method that transforms any basis function to the ones which satisfy the essential boundary condition.

In SSM, the distance function is approximated and the introduction of a weighting function increases the computational cost. Compared to IGA, SSM cannot guarantee geometric exactness because the distance function is normally approximated by the shape functions.

\subsection{Scaled Boundary Finite Element Method}
The Scaled Boundary Finite Element Method (SBFEM) of Song and Wolf \cite{songwolf1997} is a semi-analytical numerical method which does not require fundamental solutions. In SBFEM, the Cartesian coordinate system is transformed into a curvilinear coordinate system. The shape function is only dependent on the location of the surface. In two dimensional cases, the displacement is discretised by

 \begin{equation}
 {u(\xi,\eta)}=[N(\eta)]{u(\xi)}
 \end{equation}
 \noindent where $\xi$ is the radial coordinate and $\eta$ is the tangential coordinate, $N(\eta)$ is the shape function which only depends on $\eta$, hence only the boundary of the domain needs to be discretised. Substituting the above equation into the  weak form of the governing equations leads to an ordinary differential equation in $\xi$ which can be solved analytically. SBFEM retains the accuracy and inherits the advantages of the FEM and is thus suitable for non-linear problems. It also shares some of its properties with the BEM, that is, only the boundary of the domain needs to be meshed. It offers great facility for infinite domains problems and crack propagation problems. One special advantage of the SBFEM for  fracture problems is that the crack surface does not need to be meshed. However, the SBFEM is not available for  complex geometries in which cases the scaling center used to define the curvilinear coordinate system is not easily defined.

\section{Conclusions}
\label{section:future-direct-chall}

This paper critically reviewed some of the recent advances in attempting to integrate CAD and Analysis, with a special emphasis on isogeometric analysis and NURBS-enhanced finite element methods.

From this standpoint, it was said that isogeometric analysis, which shares a number of features with meshfree methods, although it has the important added advantage of conserving exact geometries,  appears like a strong method to achieve tighter CAD-Analysis integration. For this to happen, however, a few points need to be addressed:

\begin{itemize}
\item The basis functions must be the same as that used to describe the geometry, thus losing the flexibility to take advantage of other types basis functions, which has particular use in some cases, such as in the spectral element method or in structural optimisation.
\item A general and efficient algorithm  needs  to be developed to construct the parameterisation of the interior of the domain, based on surface information.
\item Mixed isogeometric methods should be studied in order to successfully deal with instabilities and construct multi-field isogeometric elements.
\item CAD uses in general multi-patch NURBS boundary descriptions, which often lead to holes, gaps, trimmed surfaces, etc. T-splines do not have this drawback, and single patch T-splines are able to provide a description for most complicated geometries seamlessly. However, the implementation of T-splines remains a challenge and it would most likely require the CAD community to embrace T-splines for T-spline based isogeometric analysis to have significant impact on current industrial practice.
\item IGA divide the elements in the parametric space, so it is difficult to capture heterogenous materials because they are scattered in the physical space. 
\end{itemize}

It seems likely that some combination of the various methods which were examined in this paper will provide winning strategies towards CAD-Analysis integration. The use of CAD-functions with Cartesian/Oct-tree meshes, for example, appears like an interesting avenue for investigation.

Our open source codes for 3D enriched isogeometric analysis, extended and smoothed extended finite element method and extended meshfree methods can be downloaded from our Source Forge Page at \url{http://cmechanicsos.users.sourceforge.net/}.

\newpage
\bibliographystyle{ECTstyle}
\bibliography{bibliography}

\end{document}